# CONDITIONS FOR RAPID MIXING OF PARALLEL AND SIMULATED TEMPERING ON MULTIMODAL DISTRIBUTIONS

By Dawn B. Woodard , Scott C. Schmidler  and Mark Huber

*Duke University*

We give conditions under which a Markov chain constructed via parallel or simulated tempering is guaranteed to be rapidly mixing, which are applicable to a wide range of multimodal distributions arising in Bayesian statistical inference and statistical mechanics. We provide lower bounds on the spectral gaps of parallel and simulated tempering. These bounds imply a single set of sufficient conditions for rapid mixing of both techniques. A direct consequence of our results is rapid mixing of parallel and simulated tempering for several normal mixture models, and for the mean-field Ising model.

**1. Introduction.** Stochastic sampling methods have become ubiquitous in statistics, computer science and statistical physics. When independent samples from a target distribution are difficult to obtain, a widely-applicable alternative is to construct a Markov chain having the target distribution as its limiting distribution. Sample paths from simulating the Markov chain then yield laws of large numbers and often central limit theorems [20, 26], and thus are widely used for Monte Carlo integration and approximate counting. Such Markov chain Monte Carlo (MCMC) methods have revolutionized computation in Bayesian statistics [9], provided significant breakthroughs in theoretical computer science [10] and become a staple of physical simulations [2, 18].

A common difficulty arising in the application of MCMC methods is that many target distributions arising in statistics and statistical physics are strongly multimodal; in such cases, the Markov chain can take an impractically long time to reach stationarity. Since the most commonly used MCMC algorithms construct reversible Markov chains, or can be made reversible without significant alteration, the convergence rate is bounded by the spectral gap of the transition operator (kernel). A variety of techniques have been developed to obtain bounds on the spectral gap of reversible Markov chains [6, 12, 24, 25]. For multimodal or other target distributions where the state space can be partitioned into high probability subsets between which the kernel rarely moves, the spectral gap will be small.

Two of the most popular and empirically successful MCMC algorithms for multimodal problems are Metropolis-coupled MCMC or *parallel tempering* [7] and *simulated tempering* [8, 17]. Adequate theoretical characterization of chains constructed in such a manner is therefore of significant interest. Toward this end, Zheng [29] bounds the spectral gap of simulated tempering below by a multiple of the spectral gap of parallel tempering, with a multiplier depending on a measure of overlap between distributions at adjacent temperatures. Madras and Piccioni [13] analyze a variant of simulated tempering as a mixture of the component chains at each temperature.

Madras and Randall [14] develop decomposition theorems for bounding the spectral gap of a Markov chain, then use those theorems to bound the mixing (equiv., convergence) of simulated tempering in terms of the slowest mixing of the tempered chains. If Metropolis–Hastings mixes slowly on the original (untempered) distribution, their bound cannot be used to show rapid mixing of simulated tempering.

However, rapid mixing of simulated tempering has been shown for several specific multimodal distributions for which local Metropolis–Hastings mixes slowly. Madras and Zheng [16] bound the spectral gap of parallel and simulated tempering on two examples, the "exponential valley" density and the mean-field Ising model. They use the decomposition theorems of [14]. However, unlike [14], they decompose the state spaces of their examples into two symmetric halves. Then they bound the mixing of parallel and simulated tempering in terms of the mixing of Metropolis–Hastings within each half. Since for these examples Metropolis–Hastings is rapidly mixing on each half of the space, their bounds show rapid mixing of parallel and simulated tempering. This is in contrast to the standard (untempered) Metropolis–Hastings chain, which is torpidly mixing. Here "torpid mixing" means that the spectral gap decreases exponentially as a function of the problem size, while "rapid mixing" means that it decreases polynomially. The torpid/rapid mixing distinction is a measure of the computational efficiency of the algorithm.

The results of [16] are extended by Bhatnagar and Randall [1] to show the rapid mixing of parallel and simulated tempering on an asymmetric version of the exponential valley density and the rapid mixing of a variant of parallel tempering on the mean-field Ising model with external field. These authors also show that parallel and simulated tempering are torpidly mixing on the mean-field Potts model with $q = 3$, regardless of the number and choice of temperatures.

We generalize the decomposition approach of [16] and [1] to obtain lower bounds on the spectral gaps of parallel and simulated tempering for any target distribution, defined on any state space, and any choice of temperatures (Theorem 3.1 and Corollary 3.1). Conceptually, we partition the state space into subsets on which the target density is unimodal. Then we bound the spectral gap of parallel and simulated tempering in terms of the mixing within each unimodal subset and the mixing among the subsets. Since Metropolis–Hastings for a unimodal distribution is often rapidly mixing, these bounds can be tighter than the simulated tempering bound of [14].

Our bounds imply a set of conditions under which parallel and simulated tempering chains are guaranteed to be rapidly mixing. The first is that Metropolis–Hastings is rapidly mixing when restricted to any one of the unimodal subsets. The challenge is then to ensure that the tempering chain is able to cross between the modes efficiently. In order to guarantee rapid mixing of the tempering chain, the second condition is that the highest-temperature chain mixes rapidly among the unimodal subsets. The third is that the overlap between distributions at adjacent temperatures decreases no more than polynomially in the problem size, which is necessary in order to mix rapidly among the temperatures. In the case where the modes are symmetric (as defined in Section 4.1), these conditions guarantee rapid mixing. We give two examples where they hold: an equally weighted mixture of normal distributions in $\mathbb{R}^M$ with identity covariance matrices (Section 4.1.2), and the mean-field Ising model (Section 4.1.1). Mixtures of normal distributions are of interest due to the fact that they closely approximate many multimodal distributions encountered in statistical practice.

When the modes are asymmetric, the three conditions above are not enough to guarantee rapid mixing, as implied by counterexamples given in the companion paper by Woodard, Schmidler and Huber [28]. In Section 4.2, we obtain an additional (fourth) condition that guarantees rapid mixing in the general case, and use this to show rapid mixing of parallel and simulated tempering for a mixture of normal distributions with unequal weights.

**2. Preliminaries.** Consider a measure space $(\mathcal{X}, \mathcal{F}, \lambda)$. Often $\mathcal{X}$ is countable and $\lambda$ is counting measure, or $\mathcal{X} = \mathbb{R}^d$ and $\lambda$ is Lebesgue measure, but more general spaces are possible. When we refer to a subset $A \subset \mathcal{X}$, we will implicitly assume that it is measurable with respect to $\lambda$. In order to draw samples from a distribution $\mu$ on $(\mathcal{X}, \mathcal{F})$, one may simulate a Markov chain that has $\mu$ as its limiting distribution, as we now describe. Let $P$ be a transition kernel on $\mathcal{X}$, defined as in [26], which operates on distributions on the left, so that for any distribution $\mu$:

$$(\mu P)(A) = \int \mu(dx) P(x, A) \qquad \forall A \subset \mathcal{X}.$$

If $\mu P = \mu$, then call $\mu$ a stationary distribution of $P$. One way of finding a transition kernel with stationary distribution $\mu$ is by constructing it to be reversible with respect to $\mu$, as we now describe.

$P$ operates on real-valued functions $f$ on the right, so that for any such $f$,

$$(Pf)(x) = \int f(y)P(x, dy) \qquad \forall x \in \mathcal{X}.$$

Define the inner product $(f, g)_\mu = \int f(x)g(x)\mu(dx)$ and denote by $L_2(\mu)$ the set of functions $f$ such that $(f, f)_\mu < \infty$. $P$ is called *reversible* with respect to $\mu$ if $(f, Pg)_\mu = (Pf, g)_\mu$ for all $f, g \in L_2(\mu)$ and *nonnegative definite* if $(Pf, f)_\mu \geq 0$ for all $f \in L_2(\mu)$. If $P$ is reversible with respect to $\mu$, then $\mu$ is easily seen to be a stationary distribution of $P$. We will primarily be interested in the case where $\mu$ has a density $\pi$ with respect to $\lambda$, and we define $\pi[A] = \mu(A)$ and define $(f, g)_\pi$, $L_2(\pi)$, and $\pi$-reversibility to be the same as for $\mu$.

If $P$ is $\phi$-irreducible and aperiodic (defined as in [22]), nonnegative definite and $\mu$-reversible, then the Markov chain with transition kernel $P$ converges in distribution to $\mu$ at a rate bounded by the *spectral gap*:

$$(1) \qquad \mathbf{Gap}(P) = \inf_{\substack{f \in L_2(\mu) \\ \mathrm{Var}_\mu(f) > 0}} \left( \frac{\mathcal{E}(f, f)}{\mathrm{Var}_\mu(f)} \right),$$

where $\mathcal{E}(f, f)$ is the *Dirichlet form* $(f, (I - P)f)_\mu$, and $\mathrm{Var}_\mu(f)$ is the variance $(f, f)_\mu - (f, 1)_\mu^2$. That is, for every distribution $\mu_0$ having a density with respect to $\mu$, there is some $h(\mu_0) > 0$ such that $\forall n \in \mathbb{N}$,

$$\begin{aligned}
(2) \qquad \|\mu_0 P^n - \mu\|_{\mathrm{TV}} &= \sup_{A \subset \mathcal{X}} |\mu_0 P^n(A) - \mu(A)| \\
&\leq h(\mu_0)[1 - \mathbf{Gap}(P)]^n \leq h(\mu_0)e^{-n\,\mathbf{Gap}(P)},
\end{aligned}$$

where $\|\cdot\|_{\mathrm{TV}}$ is the total variation distance, and the first inequality comes from functional analysis [15, 21, 23]. When $\mathbf{Gap}(P) > 0$, the chain is called geometrically ergodic (see, e.g., [22, 23]) and $\mathbf{Gap}(P)$ provides a nontrivial bound on the convergence rate. Under these conditions, we can obtain samples from $\mu$ by simulating $P$ until it has converged arbitrarily close to $\mu$. We now describe a common way of constructing a transition kernel that is reversible with respect to a particular density of interest $\pi$.

2.1. *Metropolis–Hastings.* Consider a transition kernel $P(w, dz)$ (the "proposal kernel") having a density $p(w, \cdot)$ with respect to $\lambda$ for each $w$, so that $P(w, dz) = p(w, z)\lambda(dz)$, and define the "Metropolis–Hastings transition kernel for $P$ with respect to $\pi$" as follows. If the current state is $w$, propose a move $z$ according to $P(w, \cdot)$, accept the move with probability

$$\rho(w, z) = \min\left\{1, \frac{\pi(z)p(z, w)}{\pi(w)p(w, z)}\right\}$$

and otherwise reject. Denote the resulting transition kernel by $P_{\mathrm{MH}}$; it is easily seen to be reversible with respect to $\pi$.

2.2. *Parallel and simulated tempering.* If the Metropolis–Hastings proposal kernel moves only locally in the space, and if $\pi$ has more than one mode, then the chain $P_{\mathrm{MH}}$ may move between the modes of $\pi$ infrequently. Tempering is a modification of Metropolis–Hastings where the density of interest $\pi$ is "flattened" in order to facilitate movement among the modes of $\pi$.

The parallel tempering algorithm [7] simulates parallel Markov chains defined on *tempered* densities $\pi_k(z) \propto \pi(z)^{\beta_k}$ for $z \in \mathcal{X}$, where $(\beta_0, \beta_1, \ldots, \beta_N)$ is a sequence of "inverse temperature" parameters chosen to satisfy $0 \leq \beta_0 < \cdots < \beta_N = 1$ and $\int \pi(z)^{\beta_0} \lambda(dz) < \infty$. The choice of inverse temperatures is flexible; a common choice is a geometric progression, and [19] provides an asymptotic optimality result to support this. The chains occasionally "swap" the states of adjacent temperature levels, resulting in a single Markov chain with state $x = (x_{[0]}, \ldots, x_{[N]})$ on the space $\mathcal{X}_{\mathrm{pt}} = \mathcal{X}^{N+1}$. Swaps are accepted according a Metropolis criteria which preserves the joint density

$$\pi_{\mathrm{pt}}(x) = \prod_{k=0}^{N} \pi_k(x_{[k]}), \qquad x \in \mathcal{X}_{\mathrm{pt}},$$

with product measure $\lambda_{\mathrm{pt}}(dx) = \prod_{k=0}^{N} \lambda(dx_{[k]})$. It is easy to see that the marginal density of $x_{[N]}$ under stationarity is $\pi$, the density of interest.

For concreteness, we consider the following specific interleaving of swap and update moves, and we add to each move a $1/2$ holding probability to guarantee nonnegative definiteness. The update move chooses $k$ uniformly on $\{0, \ldots, N\}$ and updates $x_{[k]}$ according to some $\pi_k$-reversible transition kernel $T_k$, yielding a transition kernel $T$ on $\mathcal{X}_{\mathrm{pt}}$:

$$T(x, dy) = \frac{1}{2(N+1)} \sum_{k=0}^{N} T_k(x_{[k]}, dy_{[k]}) \delta_{x_{[-k]}}(y_{[-k]}) \, dy_{[-k]}, \qquad x, y \in \mathcal{X}_{\mathrm{pt}},$$

where $x_{[-k]} = (x_{[0]}, x_{[1]}, \ldots, x_{[k-1]}, x_{[k+1]}, \ldots, x_{[N]})$ and $\delta$ is Dirac's delta function. Often $T_k$ is a Metropolis–Hastings kernel with respect to $\pi_k$.

The swap move $Q$ samples $k$ uniformly from $\{0, \ldots, N-1\}$ and proposes exchanging the value of $x_{[k]}$ with that of $x_{[k+1]}$. The proposed state, denoted $(k, k+1)x$, is accepted according to the Metropolis criteria preserving $\pi_{\mathrm{pt}}$:

$$\rho(x, (k, k+1)x) = \min\left\{1, \frac{\pi_k(x_{[k+1]})\pi_{k+1}(x_{[k]})}{\pi_k(x_{[k]})\pi_{k+1}(x_{[k+1]})}\right\}$$

so that for any $A \subset \mathcal{X}_{\text{pt}}$ and $x \in \mathcal{X}_{\text{pt}}$,

$$Q(x, A) = \frac{1}{2N} \sum_{k=0}^{N-1} \mathbf{1}_A((k, k+1)x) \rho(x, (k, k+1)x)$$
$$+ \mathbf{1}_A(x) \left[ 1 - \frac{1}{2N} \sum_{k=0}^{N-1} \rho(x, (k, k+1)x) \right],$$

where $\mathbf{1}_A$ is the indicator function of the set $A$. Both $T$ and $Q$ are easily seen to be reversible with respect to $\pi_{\text{pt}}$ by construction, and nonnegative definite due to their $1/2$ holding probability. We define the parallel tempering chain $P_{\text{pt}} = QTQ$ which performs two swapping moves for each update move. It is easily verified that $P_{\text{pt}}$ is nonnegative definite and reversible with respect to $\pi_{\text{pt}}$, so the convergence of the parallel tempering chain to $\pi_{\text{pt}}$ may be bounded using the spectral gap of $P_{\text{pt}}$.

Note that the definitions of $T$ and $Q$ do not rely on $\pi_k \propto \pi^{\beta_k}$, and indeed we may specify distributions $\pi_k$ in any convenient way subject to $\pi_N = \pi$. We refer to the resulting chain as a *swapping chain*, and denote the corresponding state space, measure, transition kernel and associated stationary density by $\mathcal{X}_{\text{sc}}$, $\lambda_{\text{sc}}$, $P_{\text{sc}}$ and $\pi_{\text{sc}}$, respectively. Although the terms "parallel tempering chain" and "swapping chain" are used interchangeably in the computer science literature, we follow the statistics and physics literature in defining a parallel tempering chain as using tempered distributions, and define a "swapping chain" as using arbitrary distributions.

The related technique of simulated tempering [8, 17] has state $(z, k) \in \mathcal{X}_{\text{st}} = \mathcal{X} \otimes \{0, \ldots, N\}$ and stationary density

$$\pi_{\text{st}}(z, k) = \frac{1}{N+1} \pi_k(z), \qquad (z, k) \in \mathcal{X}_{\text{st}},$$

with two move types: the first $(T')$ updates $z \in \mathcal{X}$ according to $T_k$, conditional on $k$, and the second $(Q')$ samples the level $k$ from its conditional distribution given $z$. Once again, a holding probability of $1/2$ is added to both $T'$ and $Q'$. The transition kernel is then specified as $P_{\text{st}} = Q'T'Q'$. For a lack of separate terms, we use "simulated tempering" to mean any such chain $P_{\text{st}}$, regardless of whether or not the densities $\pi_k$ are tempered versions of $\pi$.

**3. Lower bounds on the spectral gaps of swapping and simulated tempering chains.** Two key results of the current paper are lower bounds on the spectral gaps of swapping and simulated tempering chains. These bounds imply the conditions for rapid mixing given in Section 4. The bounds are in terms of several quantities. Informally, the first quantity measures how well each chain $T_k$ mixes when restricted to each unimodal subset. The second is how well the highest-temperature chain $T_0$ mixes among the subsets. The third is the overlap of the distributions of adjacent levels, and the fourth concerns the probability of each unimodal subset as a function of the (inverse) temperature.

In order to bound the spectral gap of a swapping or simulated tempering chain in terms of the mixing of the chain within each subset and the mixing of the chain among the subsets, we use a state space decomposition result due to Caracciolo, Pelissetto and Sokal [3] and first published by Madras and Randall [14]. As in [14], we will use the following definitions.

For any transition kernel $P$ reversible with respect to a distribution $\mu$ and any subset $A$ of the state space of $P$, define the restriction of $P$ to $A$ as

$$(3) \qquad P|_A(x, B) = P(x, B) + \mathbf{1}_B(x)P(x, A^c) \qquad \text{for } x \in A, B \subset A.$$

Note that $P|_A$ is reversible with respect to $\mu|_A$, the restriction of $\mu$ to $A$. Now take any partition $\mathcal{A} = \{A_j : j = 1, \ldots, J\}$ of the state space of $P$ such that $\mu(A_j) > 0$ for all $j$, and define the projection matrix of $P$ with respect to $\mathcal{A}$ as

$$(4) \qquad \bar{P}(i, j) = \frac{1}{\mu(A_i)} \int_{A_i} \int_{A_j} P(x, dy)\mu(dx), \qquad i, j \in \{1, \ldots, J\}.$$

Note that $\bar{P}$ is reversible with respect to the distribution on $j \in \{1, \ldots, J\}$ taking value $\mu(A_j)$, and irreducible if $P$ is.

Now consider a swapping or simulated tempering chain defined as in Section 2.2 for some density of interest $\pi$ on a measure space $(\mathcal{X}, \mathcal{F}, \lambda)$, with $\pi_k$-reversible transition kernels $T_k$. Let $\mathcal{A}$ be any partition of $\mathcal{X}$ such that $\pi_k[A_j] > 0$ for all $k$ and $j$. The first quantity in our bound is the minimum over $k$ and $j$ of $\mathbf{Gap}(T_k|_{A_j})$, which measures how well each chain $T_k$ mixes within each partition element. The partition would typically be chosen so that this quantity is large; in our examples, we choose $\mathcal{A}$ so that $\pi|_{A_j}$ is unimodal (contains a single local mode) for each $j$.

Next, we consider how well the chain $T_0$ mixes among the partition elements. Let $\bar{T}_0$ be the projection matrix of $T_0$ with respect to $\mathcal{A}$; the second quantity in our bound is $\mathbf{Gap}(\bar{T}_0)$. Since $\mathcal{A}$ is finite, this is one minus the second-largest eigenvalue of $\bar{T}_0$.

The third quantity is the *overlap* of $\{\pi_k : k = 0, \ldots, N\}$ with respect to $\mathcal{A}$, defined as

$$(5) \qquad \delta(\mathcal{A}) = \min_{\substack{|k-l|=1 \\ j \in \{1, \ldots, J\}}} \left[ \int_{A_j} \min\{\pi_k(z), \pi_l(z)\}\lambda(dz) \right] \Big/ \pi_k[A_j].$$

The quantity $\delta(\mathcal{A})$ controls the rate of temperature changes in simulated tempering. For the swapping chain, note that for any $i, j \in \{1, \ldots, J\}$ and any $k \in \{0, \ldots, N-1\}$, the marginal probability at stationarity of accepting a proposed swap between $x_{[k]} \in A_i$ and $x_{[k+1]} \in A_j$ is

$$(6) \qquad \frac{\int_{z \in A_i} \int_{w \in A_j} \min\{\pi_k(z)\pi_{k+1}(w), \pi_k(w)\pi_{k+1}(z)\}\lambda(dw)\lambda(dz)}{\pi_k[A_i]\pi_{k+1}[A_j]} \geq \delta(\mathcal{A})^2.$$

We will show that our overlap quantity $\delta(\mathcal{A})$ is bounded below by the overlap used in Madras and Randall [14] and Zheng [29], and that our definition is equal to theirs in the case of $\pi$ symmetric (as defined in Section 4.1).

The fourth and final quantity concerns the probability of a single partition element under $\pi_k$, as a function of $k$, for each partition element:

$$(7) \qquad \gamma(\mathcal{A}) = \min_{j \in \{1, \ldots, J\}} \prod_{k=1}^{N} \min\left\{1, \frac{\pi_{k-1}[A_j]}{\pi_k[A_j]}\right\}.$$

Note that for any $j \in \{1, \ldots, J\}$ and any $k, l \in \{0, \ldots, N\}$ such that $k < l$, $\pi_k[A_j] \geq \gamma(\mathcal{A})\pi_l[A_j]$. If $\pi_k[A_j]$ is a monotonic function of $k$ for each $j$ (which need not hold for tempered distributions), then $\gamma(\mathcal{A})$ simplifies to

$$\min_{j \in \{1, \ldots, J\}} \frac{\pi_0[A_j]}{\pi_N[A_j]}.$$

With these definitions, the following theorem bounds the spectral gap of the swapping chain.

THEOREM 3.1. *Given any partition* $\mathcal{A} = \{A_j : j = 1, \ldots, J\}$ *of* $\mathcal{X}$ *such that* $\pi_k[A_j] > 0$ *for all* $k$ *and* $j$, *and given* $\delta(\mathcal{A})$ *as in* (5) *and* $\gamma(\mathcal{A})$ *as in* (7),

$$\mathbf{Gap}(P_{sc}) \geq \left(\frac{\gamma(\mathcal{A})^{J+3}\delta(\mathcal{A})^2}{2^{12}(N+1)^4 J^3}\right)\mathbf{Gap}(\bar{T}_0) \min_{k,j} \mathbf{Gap}(T_k|_{A_j}).$$

In particular, the bound holds for parallel tempering with $\pi_k \propto \pi^{\beta_k}$. Theorem 3.1 will be proven in Section 6. Note that

$$
\begin{aligned}
\delta(\mathcal{A}) &= \min_{\substack{k \in \{0, \ldots, N-1\} \\ j \in \{1, \ldots, J\}}} \frac{\int_{A_j} \min\{\pi_k(z), \pi_{k+1}(z)\}\lambda(dz)}{\max\{\pi_k[A_j], \pi_{k+1}[A_j]\}} \\
(8) \qquad &\leq \min_{k \in \{0, \ldots, N-1\}} \frac{\sum_j \int_{A_j} \min\{\pi_k(z), \pi_{k+1}(z)\}\lambda(dz)}{\max\{\sum_j \pi_k[A_j], \sum_j \pi_{k+1}[A_j]\}} \\
&= \min_{k \in \{0, \ldots, N-1\}} \int \min\{\pi_k(z), \pi_{k+1}(z)\}\lambda(dz).
\end{aligned}
$$

The final expression for $\delta(\mathcal{A})$ is the definition of overlap that is used in [14] and [29]. Therefore, Theorem 3 of [29], along with our Theorem 3.1, implies the following bound for simulated tempering:

COROLLARY 3.1. *Let $P_{st}$ be the simulated tempering chain defined with the same $N$ and same set of densities $\pi_k$ as the swapping chain. Then*

$$\mathbf{Gap}(P_{st}) \geq \left(\frac{\gamma(\mathcal{A})^{J+3}\delta(\mathcal{A})^3}{2^{14}(N+1)^5 J^3}\right)\mathbf{Gap}(\bar{T}_0)\min_{k,j}\mathbf{Gap}(T_k|_{A_j}).$$

**4. Examples of rapid mixing.** We will show rapid mixing of parallel and simulated tempering for several examples by applying Theorem 3.1 and Corollary 3.1. We are particularly interested in cases where Metropolis–Hastings with local proposals is torpidly mixing due to the multimodality of the target density $\pi$. To show rapid mixing of tempering for a case where the number of modes $J$ is fixed, we choose a set of temperatures the number of which grows at most polynomially in the problem size. We then show that each chain $T_k$ is rapidly mixing when restricted to each unimodal subset [meaning that $\min_{j,k}\mathbf{Gap}(T_k|_{A_j})$ decreases at most polynomially in the problem size], and that the highest-temperature chain mixes rapidly among the subsets. Additionally, we show that the overlap of the parallel or simulated tempering chain decreases at most polynomially in the problem size, which is needed in order to mix rapidly among the temperatures. In the case where $\pi$ is symmetric with respect to $\mathcal{A}$ (to be defined), we will see that $\gamma(\mathcal{A}) = 1$, so the above conditions imply that parallel tempering is rapidly mixing (by Theorem 3.1). The same conditions imply the rapid mixing of simulated tempering (by Corollary 3.1) using the same tempered densities.

The above conditions are necessary as well as sufficient for rapid mixing in the symmetric case. Simple examples can be constructed to show the necessity of each condition; Woodard, Schmidler and Huber [28] give an example of a symmetric distribution for which the first two conditions hold, but the condition on the overlap fails, and parallel and simulated tempering are torpidly mixing.

In the case of general (not necessarily symmetric) $\pi$, the above conditions are insufficient to guarantee rapid mixing; counterexamples are given in [28]. In the general case, we must additionally show that $\gamma(\mathcal{A})$ decreases at most polynomially in the problem size.

4.1. *Examples of rapid mixing on symmetric distributions.* Recall that the target density $\pi$ is defined on a state space $\mathcal{X}$ with measure $\lambda$. Define $\pi$ to be *symmetric* with respect to a partition $\{A_j : j = 1, \ldots, J\}$ of $\mathcal{X}$ if for every pair of partition elements $A_i$, $A_j$ there is some $\lambda$-measure-preserving bijection $f_{ij}$ from $A_i$ to $A_j$ that preserves $\pi$. Note that when $\pi$ is symmetric with respect to $\{A_j : j = 1, \ldots, J\}$, the inequality in (8) is an equality. Additionally, $\pi_k[A_j] = 1/J$ for all $k \in \{0, \ldots, N\}$ and $j \in \{1, \ldots, J\}$, so $\gamma(\mathcal{A}) = 1$. We will give examples of symmetric $\pi$ for which parallel and simulated tempering are rapidly mixing.

4.1.1. *The mean field Ising model.* For each $M \in \mathbb{N}$, the mean field Ising model is defined for $z \in \mathcal{X} = \{-1, +1\}^M$ as:

$$
(9) \qquad \pi(z) = \frac{1}{Z} \exp\left\{ \frac{\alpha}{2M} \left( \sum_{i=1}^{M} z_i \right)^2 \right\},
$$

where $Z = \sum_z \exp\{\alpha(\sum_i z_i)^2/(2M)\}$. The single-site proposal kernel $S$ chooses $i \in \{1, \dots, M\}$ uniformly at random and proposes switching the sign of $z_i$. Metropolis–Hastings for $S$ with respect to $\pi$ is torpidly mixing for $\alpha > 1$, as is straightforward to show using conductance (as defined in [12, 25]).

Taking $N = M$, $\beta_k = k/N$ and $T_k$ equal to Metropolis–Hastings for $S$ with respect to $\pi_k \propto \pi^{\beta_k}$, it is shown by Madras and Zheng [16] that parallel and simulated tempering are rapidly mixing. We will show that this is also a consequence of our Theorem 3.1 and Corollary 3.1.

As in [16], partition $\mathcal{X}$ into $A_1 = \{z \in \mathcal{X} : \sum_i z_i < 0\}$ and $A_2 = \{z \in \mathcal{X} : \sum_i z_i \geq 0\}$. Restricting to $M$ odd, the density $\pi$ is clearly symmetric with respect to the partition.

Using the fact that $\pi_0$ is uniform, it is straightforward to show that $T_0$ is rapidly mixing. Since $T_0$ is rapidly mixing, so is $\bar{T}_0$ (see Theorem 5.2). Additionally, it is shown in [16] that the minimum over $k$ and $j$ of $\mathbf{Gap}(T_k|_{A_j})$ is polynomially decreasing in $M$.

Note that for any $z \in \mathcal{X}$ and any $k \in \{0, \dots, N-1\}$,

$$
\pi(z)^{\beta_{k+1} - \beta_k} = \pi(z)^{1/M} \in \left[ \frac{1}{Z^{1/M}}, \frac{1}{Z^{1/M}} \exp\{\alpha/2\} \right].
$$

Therefore,

$$
\frac{\pi_{k+1}(z)}{\pi_k(z)} = \pi(z)^{\beta_{k+1} - \beta_k} \left( \frac{\sum_{w \in \mathcal{X}} \pi(w)^{\beta_k}}{\sum_{w \in \mathcal{X}} \pi(w)^{\beta_{k+1}}} \right) \in [\exp\{-\alpha/2\}, \exp\{\alpha/2\}]
$$

which implies that

$$
\sum_{z \in \mathcal{X}} \min\{\pi_k(z), \pi_{k+1}(z)\} = \sum_{z \in \mathcal{X}} \pi_k(z) \min\left\{ 1, \frac{\pi_{k+1}(z)}{\pi_k(z)} \right\}
$$
$$
\geq \exp\{-\alpha/2\}.
$$

Recalling that for $\pi$ symmetric, the inequality in (8) is an equality, $\delta(\mathcal{A})$ is bounded below by a constant for all $M$. Therefore, by Theorem 3.1 and Corollary 3.1, the parallel or simulated tempering chain is rapidly mixing.

4.1.2. *A symmetric normal mixture.* Many multimodal distributions arising in statistics are well approximated by mixtures of normal distributions. We will analyze the mixing of parallel and simulated tempering on several two-component mixtures of normal distributions in $\mathbb{R}^M$. For any length-$M$ vector $\nu$, $M \times M$ covariance matrix $\Sigma$, and $z \in \mathbb{R}^M$, let $N_M(z; \nu, \Sigma)$ be the density of a multivariate normal distribution in $\mathbb{R}^M$ with mean $\nu$ and covariance $\Sigma$, evaluated at $z$. Let $1_M$ denote the vector of $M$ ones, and $I_M$ denote the $M \times M$ identity matrix. Take any $b > 0$ and any sequence $a_1, a_2, \ldots$ such that $a_M \in (0, 1)$ for each $M$, and consider the following weighted mixture of two normal densities in $\mathbb{R}^M$:

$$(10) \qquad \pi(z) = a_M N_M(z; -b1_M, I_M) + (1 - a_M) N_M(z; b1_M, I_M).$$

Let $S$ be the proposal kernel that is uniform on the ball of radius $M^{-1}$ centered at the current state. Partition $\mathcal{X}$ into $A_1 = \{z : \sum_i z_i < 0\}$ and $A_2 = \{z : \sum_i z_i \geq 0\}$. For technical reasons, we will use the following approximation to $\pi$:

$$(11) \qquad \tilde{\pi}(z) \propto a_M N_M(z; -b1_M, I_M) \mathbf{1}_{A_1}(z) + (1 - a_M) N_M(z; b1_M, I_M) \mathbf{1}_{A_2}(z)$$

which truncates the overlapping portions of the tails. This simplification does not alter the mixing properties of Metropolis–Hastings or tempering, with the exception of the case where either $a_M$ or $1 - a_M$ approaches zero very quickly so that $\pi$ is unimodal for large $M$ while $\tilde{\pi}$ is bimodal. However, we are interested only in the situation where $\pi$ is bimodal, causing poor mixing of Metropolis–Hastings and suggesting the more efficient application of tempering; for this purpose $\tilde{\pi}$ is equivalent to $\pi$.

Metropolis–Hastings for $S$ with respect to the density

$$\tilde{\pi}|_{A_1}(z) \propto N_M(z; -b1_M, I_M) \mathbf{1}_{A_1}(z)$$

or with respect to

$$\tilde{\pi}|_{A_2}(z) \propto N_M(z; b1_M, I_M) \mathbf{1}_{A_2}(z)$$

is rapidly mixing in $M$, as implied by results in Kannan and Li [11] (details are given in [27]). However, we will show that Metropolis–Hastings for $S$ with respect to $\tilde{\pi}$ is torpidly mixing. Then we will specify a set of inverse temperatures which yield rapidly mixing parallel and simulated tempering chains for $\tilde{\pi}$.

Consider Metropolis–Hastings for $S$ with respect to $\tilde{\pi}$. Note that the boundary of $A_1$ with respect to the Metropolis–Hastings kernel is the set of $z \in A_1$ within distance $M^{-1}$ of the set $A_2$. It is straightforward to show that the probability of this boundary under $\tilde{\pi}|_{A_1}$ decreases exponentially as a function of $M$. Similarly, the probability of the boundary of $A_2$ under $\tilde{\pi}|_{A_2}$ decreases exponentially in $M$. Therefore, the conductance of the set $A_1$ (where conductance is defined as in [12, 25]) is exponentially decreasing, which implies that Metropolis–Hastings is torpidly mixing.

For any $\beta$, define the tempered density $\tilde{\pi}_\beta \propto \tilde{\pi}^\beta$. Note that for any $\beta$,

$$\tilde{\pi}_\beta(z) \propto a_M^\beta N_M(z; -b1_M, \beta^{-1} I_M) \mathbf{1}_{A_1}(z)$$
$$\qquad + (1 - a_M)^\beta N_M(z; b1_M, \beta^{-1} I_M) \mathbf{1}_{A_2}(z).$$

The normalizing constant is $[a_M^\beta + (1-a_M)^\beta]\Phi(b\sqrt{M}\beta^{1/2})$, where $\Phi$ is the cumulative normal distribution function in one dimension.

Metropolis–Hastings for $S$ with respect to $\tilde{\pi}_{M^{-1}}$ mixes rapidly between $A_1$ and $A_2$, shown as follows. The probability of the boundary of $A_1$ under $\tilde{\pi}_{M^{-1}}|_{A_1}$ is

$$[\Phi(b) - \Phi(b(1-M^{-1}))]/\Phi(b) \tag{12}$$

which is polynomially decreasing in $M$. The probability of the boundary of $A_2$ under $\tilde{\pi}_{M^{-1}}|_{A_2}$ is also equal to (12). It can also be shown that the marginal probability of accepting proposed moves between $A_1$ and $A_2$ is polynomially decreasing in $M$, proving rapid mixing between $A_1$ and $A_2$; the details are given in [27].

The infimum over $\beta \geq M^{-1}$ of the spectral gap of Metropolis–Hastings for $S$ with respect to $\tilde{\pi}_\beta|_{A_j}$ decreases at most polynomially in $M$ for $j = 1, 2$. This is because $\tilde{\pi}_\beta|_{A_j}$ is a normal density restricted to a convex set; a bound on the spectral gap of Metropolis–Hastings for $S$ with respect to such a density is given in [11], and this bound is polynomially decreasing in $M$. More details are given in [27].

When $a_M = 1 - a_M$, $\tilde{\pi}$ is symmetric with respect to the partition $\{A_1, A_2\}$. In this case, set $N = M$ and $\beta_k = M^{-(M-k)/M}$ (a geometric progression), and let $T_k$ be the Metropolis–Hastings kernel for $S$ with respect to $\pi_k \propto \tilde{\pi}^{\beta_k}$. With these specifications, parallel and simulated tempering are rapidly mixing, as we will show. We have already seen that $\min_{j,k} \mathbf{Gap}(T_k|_{A_j})$ is polynomially decreasing in $M$ and that $\bar{T}_0$ is rapidly mixing. Next, we will show that $\delta(\mathcal{A})$ is also polynomially decreasing in $M$.

Let $\lambda$ be Lebesgue measure in $\mathbb{R}^M$, and take any $k \in \{0, \ldots, N-1\}$. Noting that $\beta_k/\beta_{k+1} = M^{-1/M}$, we have

$$\int_{\mathcal{X}} \min\{\pi_k(z), \pi_{k+1}(z)\}\lambda(dz)$$

$$= 2\int_{A_2} \min\{\pi_k(z), \pi_{k+1}(z)\}\lambda(dz)$$

$$= \int_{A_2} \min\left\{\frac{N_M(z; b1_M, \beta_k^{-1}I_M)}{\Phi(b\sqrt{M}\beta_k^{1/2})}, \frac{N_M(z; b1_M, \beta_{k+1}^{-1}I_M)}{\Phi(b\sqrt{M}\beta_{k+1}^{1/2})}\right\}\lambda(dz)$$

$$\geq \int_{A_2} \min\{N_M(z; b1_M, \beta_k^{-1}I_M), N_M(z; b1_M, \beta_{k+1}^{-1}I_M)\}\lambda(dz)$$

$$= (2\pi)^{-M/2} \tag{13}$$

$$\times \int_{A_2} \beta_{k+1}^{M/2} \min\left\{\left(\frac{\beta_k}{\beta_{k+1}}\right)^{M/2} \exp\left\{-\frac{\beta_k}{2}\sum_i (z_i - b)^2\right\},\right.$$

$$\exp\left\{-\frac{\beta_{k+1}}{2}\sum_i(z_i-b)^2\right\}\right\}\lambda(dz)$$

$$\geq\frac{1}{\sqrt{M}}(2\pi)^{-M/2}\int_{A_2}\beta_{k+1}^{M/2}\exp\left\{-\frac{\beta_{k+1}}{2}\sum_i(z_i-b)^2\right\}\lambda(dz)$$

$$=\frac{1}{\sqrt{M}}\Phi(b\sqrt{M}\beta_{k+1}^{1/2})\geq\frac{1}{2\sqrt{M}}.$$

Therefore, $\delta(\mathcal{A})$ decreases at most polynomially in $M$. By Theorem 3.1 and Corollary 3.1, parallel and simulated tempering with this $N$ and this set of inverse temperatures are rapidly mixing in the case where $a_M=1-a_M$. This is shown in [27] for $\pi$ as well as $\tilde{\pi}$. In the case where there are some $M$ with $a_M\neq1-a_M$, we must also verify that $\gamma(\mathcal{A})$ is polynomially decreasing in $M$; this is done in the next section.

### 4.2. *Examples of rapid mixing on general distributions.* 
We now consider the case of general (not necessarily symmetric) $\pi$.

#### 4.2.1. *A weighted normal mixture.* 
Recall $\tilde{\pi}$ from (11), and assume without loss of generality that $a_M\geq1/2$. For technical reasons, we restrict $a_M/(1-a_M)$ to be exponentially bounded-above, meaning that there exists a constant $c$ such that $a_M/(1-a_M)\leq c^M$ for all $M$.

Consider the inverse temperature specification that we used for the symmetric case, with $N=M$ and the set of inverse temperatures $\{M^{-(M-k)/M}:k=0,\ldots,M\}$. Also recall the inverse temperature specification for the mean field Ising model, with $N=M$ and the set of inverse temperatures $\{k/M:k=0,\ldots,M\}$. For the mixture of normals with unequal weights, we will need both: take the set of inverse temperatures $\{M^{-(M-k)/M}:k=0,\ldots,M\}\cup\{k/M:k=1,\ldots,M\}$, so $N=2M$.

The arguments from Section 4.1.2 show that $\min_{j,k}\mathbf{Gap}(T_k|_{A_j})$ is polynomially decreasing in $M$ and that $\bar{T}_0$ is rapidly mixing. Next, we will show that $\gamma(\mathcal{A})$ and $\delta(\mathcal{A})$ are also polynomially decreasing in $M$.

Note that for any $\beta$,

$$\tilde{\pi}_\beta[A_1]=\frac{a_M^\beta}{a_M^\beta+(1-a_M)^\beta}$$

which is an increasing function of $\beta$ since $a_M\geq1/2$. Therefore, $\pi_k[A_1]$ is an increasing function of $k$ and $\pi_k[A_2]$ is a decreasing function of $k$, so

$$\gamma(\mathcal{A})=\frac{\pi_0[A_1]}{\pi_N[A_1]}\geq\frac{1}{2\pi_N[A_1]}\geq\frac{1}{2}$$

which does not depend on $M$. Also note that for any $k \in \{0, \ldots, N-1\}$,

$$\frac{\pi_{k+1}[A_2]}{\pi_k[A_2]} \geq \frac{\pi_k[A_1]\pi_{k+1}[A_2]}{\pi_{k+1}[A_1]\pi_k[A_2]} = \left(\frac{1-a_M}{a_M}\right)^{\beta_{k+1}-\beta_k}$$

$$\geq \left(\frac{1-a_M}{a_M}\right)^{1/M} \geq c^{-1}.$$

Therefore,

$$\frac{\int_{A_2} \min\{\pi_k(z), \pi_{k+1}(z)\}\lambda(dz)}{\max\{\pi_k[A_2], \pi_{k+1}[A_2]\}}$$

$$= \int_{A_2} \min\left\{\pi_k[A_2]\frac{N_M(z; b1_M, \beta_k^{-1}I_M)}{\Phi(b\sqrt{M}\beta_k^{1/2})},\right.$$

$$\left. \pi_{k+1}[A_2]\frac{N_M(z; b1_M, \beta_{k+1}^{-1}I_M)}{\Phi(b\sqrt{M}\beta_{k+1}^{1/2})}\right\}\lambda(dz)$$

$$\bigg/ \max\{\pi_k[A_2], \pi_{k+1}[A_2]\}$$

$$\geq \frac{\pi_{k+1}[A_2]\int_{A_2} \min\{N_M(z; b1_M, \beta_k^{-1}I_M), N_M(z; b1_M, \beta_{k+1}^{-1}I_M)\}\lambda(dz)}{\pi_k[A_2]}$$

$$\geq c^{-1}\int_{A_2} \min\{N_M(z; b1_M, \beta_k^{-1}I_M), N_M(z; b1_M, \beta_{k+1}^{-1}I_M)\}\lambda(dz)$$

$$\geq \frac{1}{2\sqrt{M}}c^{-1},$$

where the last inequality is from (13), using the fact that $\beta_k/\beta_{k+1} \geq M^{-1/M}$. This result, repeated for $A_1$, shows that $\delta(\mathcal{A})$ decreases at most polynomially in $M$. By Theorem 3.1, parallel and simulated tempering with this $N$ and this set of inverse temperatures are rapidly mixing on the weighted mixture of normals $\tilde{\pi}$.

**5. Tools for bounding spectral gaps.** In Sections 5.1 and 5.2, we give some results from the literature and slight extensions thereon. These results will be used in Section 6 for the proof of Theorem 3.1.

5.1. *A bound for finite state space Markov chains.* We first consider a method for finite state space Markov chains. Let $P$ and $Q$ be Markov chain transition matrices on state space $\mathcal{X}$ with $|\mathcal{X}| < \infty$, reversible with respect to densities $\pi_P$ and $\pi_Q$, respectively. Denote by $\mathcal{E}_P$ and $\mathcal{E}_Q$ the Dirichlet forms of $P$ and $Q$, and let $E_P = \{(x, y) : \pi_P(x)P(x, y) > 0\}$ and $E_Q = \{(x, y) : \pi_Q(x)Q(x, y) > 0\}$ be the edge sets of $P$ and $Q$, respectively.

For each pair $x \neq y$ such that $(x, y) \in E_Q$, fix a path $\gamma_{xy} = (x = x_0, x_1, x_2, \ldots, x_k = y)$ of length $|\gamma_{xy}| = k$ such that $(x_i, x_{i+1}) \in E_P$ for $i \in \{0, \ldots, k-1\}$. Define

$$c = \max_{(z,w) \in E_P} \left\{ \frac{1}{\pi_P(z) P(z, w)} \sum_{\gamma_{xy} \ni (z,w)} |\gamma_{xy}| \pi_Q(x) Q(x, y) \right\}.$$

THEOREM 5.1 (Diaconis and Saloff-Coste [4]).

$$\mathcal{E}_Q \leq c \mathcal{E}_P.$$

5.2. *Bounds for general state space Markov chains.* The following results hold for general state space transition kernels $P$ and $Q$, reversible with respect to distributions $\mu_P$ and $\mu_Q$ on a space $\mathcal{X}$ with countably generated $\sigma$-algebra.

THEOREM 5.2. *Let $\{A_j : j = 1, \ldots, J\}$ be any partition of $\mathcal{X}$ such that $\mu_P(A_j) > 0$ for all $j$. Define $P|_{A_j}$ as in (3) and $\bar{P}$ as in (4). For $P$ nonnegative definite,*

$$\frac{1}{2} \mathbf{Gap}(\bar{P}) \min_{j=1,\ldots,J} \mathbf{Gap}(P|_{A_j}) \leq \mathbf{Gap}(P) \leq \mathbf{Gap}(\bar{P}).$$

The bounds are a direct consequence of results published in Madras and Randall [14], as described in the Appendix of the current paper.

THEOREM 5.3 [Diaconis and Saloff-Coste (1996)]. *Take any $N \in \mathbb{N}$ and let $P_k$, $k = 0, \ldots, N$, be $\mu_k$-reversible transition kernels on state spaces $\mathcal{X}_k$. Let $P$ be the transition kernel on $\mathcal{X} = \prod_k \mathcal{X}_k$ given by*

$$P(x, dy) = \sum_{k=0}^{N} b_k P_k(x_{[k]}, dy_{[k]}) \delta_{x_{[-k]}}(y_{[-k]}) \, dy_{[-k]}, \qquad x, y \in \mathcal{X},$$

*for some set of $b_k > 0$ such that $\sum_k b_k = 1$ and where $\delta$ is Dirac's delta function. $P$ is called a product chain. It is reversible with respect to $\mu(dx) = \prod_k \mu_k(dx_{[k]})$, and*

$$\mathbf{Gap}(P) = \min_{k=0,\ldots,N} b_k \mathbf{Gap}(P_k).$$

Lemma 3.2 of [5] states this result for finite state spaces; however, the proof of that lemma holds in the general case.

LEMMA 5.1. *Let $\mu_P = \mu_Q$. If $Q(x, A \setminus \{x\}) \leq P(x, A \setminus \{x\})$ for every $x \in \mathcal{X}$ and every $A \subset \mathcal{X}$, then $\mathbf{Gap}(Q) \leq \mathbf{Gap}(P)$.*

PROOF.   As in [14], write $\mathbf{Gap}(P)$ in the form

$$\mathbf{Gap}(P) = \inf_{\substack{f \in L_2(\mu_P) \\ \mathrm{Var}_{\mu_P}(f) > 0}} \left( \frac{\iint |f(x) - f(y)|^2 \mu_P(dx) P(x, dy)}{\iint |f(x) - f(y)|^2 \mu_P(dx) \mu_P(dy)} \right)$$

and write $\mathbf{Gap}(Q)$ analogously. The result then follows immediately.   □

## 6. Proof of Theorem 3.1.

6.1. *Overview of the proof.*   As in Madras and Zheng [16], consider the space $\Sigma = \mathbb{Z}_J^{N+1}$ of possible assignments of levels to partition elements. For $x = (x_{[0]}, \dots, x_{[N]}) \in \mathcal{X}_{\mathrm{sc}}$, let the signature $s(x)$ be the vector $(\sigma_0, \dots, \sigma_N) \in \Sigma$ with

$$\sigma_k = j \qquad \text{if } x_{[k]} \in A_j \ (0 \le k \le N)$$

and for $\sigma \in \Sigma$, define

$$\mathcal{X}_\sigma = \{x \in \mathcal{X}_{\mathrm{sc}} : s(x) = \sigma\}$$

so $s$ induces a partition of $\mathcal{X}_{\mathrm{sc}}$. Define $P_\sigma = P_{\mathrm{sc}}|_{\mathcal{X}_\sigma}$, and let $\bar{P}_{\mathrm{sc}}$ be the projection matrix of $P_{\mathrm{sc}}$ with respect to the partition $\{\mathcal{X}_\sigma\}$. Since $P_{\mathrm{sc}}$ is nonnegative definite, Theorem 5.2 gives

$$(14) \qquad \mathbf{Gap}(P_{\mathrm{sc}}) \ge \frac{1}{2} \mathbf{Gap}(\bar{P}_{\mathrm{sc}}) \min_{\sigma \in \Sigma} \mathbf{Gap}(P_\sigma).$$

Theorem 3.1 then follows by deriving bounds on $\mathbf{Gap}(\bar{P}_{\mathrm{sc}})$ and $\mathbf{Gap}(P_\sigma)$.

6.2. *Bounding the spectral gap of $P_\sigma$.*   For $\sigma \in \Sigma$, consider the mixing of $P_{\mathrm{sc}}$ when restricted to the set $\mathcal{X}_\sigma$. If each of the chains $T_k$ mixes well when restricted to the set $A_{\sigma_k}$, then the product chain $T$, and thus $P_{\mathrm{sc}}$ will also mix well when restricted to $\mathcal{X}_\sigma$. Let $T_\sigma = T|_{\mathcal{X}_\sigma}$ and note that for any $x, y \in \mathcal{X}_\sigma$ with $x \ne y$,

$$T_\sigma(x, dy) = \frac{1}{2(N+1)} \sum_{k=0}^N T_k|_{A_{\sigma_k}}(x_{[k]}, dy_{[k]}) \delta_{x_{[-k]}}(y_{[-k]}) \, dy_{[-k]}.$$

Therefore, $T_\sigma$ is a product chain, and Theorem 5.3 provides its spectral gap:

$$\mathbf{Gap}(T_\sigma) = \frac{1}{2(N+1)} \min_{k \in \{0, \dots, N\}} \mathbf{Gap}(T_k|_{A_{\sigma_k}})$$

$$\ge \frac{1}{2(N+1)} \min_{\substack{k \in \{0, \dots, N\} \\ j \in \{1, \dots, J\}}} \mathbf{Gap}(T_k|_{A_j}).$$

Note that since $P_{\mathrm{sc}} = QTQ$, and since $Q$ has a $1/2$ holding probability,

$$P_\sigma(x, dy) \geq \tfrac{1}{4} T_\sigma(x, dy) \qquad \forall x, y \in \mathcal{X}_\sigma.$$

Using Lemma 5.1 we have $\mathbf{Gap}(P_\sigma) \geq \mathbf{Gap}(T_\sigma)/4$. Therefore

$$(15) \qquad \mathbf{Gap}(P_\sigma) \geq \frac{1}{8(N+1)} \min_{\substack{k \in \{0, \ldots, N\} \\ j \in \{1, \ldots, J\}}} \mathbf{Gap}(T_k|_{A_j}).$$

6.3. *Bounding the spectral gap of $\bar{P}_{\mathrm{sc}}$.* First note that $\bar{P}_{\mathrm{sc}}$ is reversible with respect to the probability mass function

$$\pi^*(\sigma) \stackrel{\mathrm{def}}{=} \pi_{\mathrm{sc}}[\mathcal{X}_\sigma] = \prod_{k=0}^N \pi_k[A_{\sigma_k}] \qquad \forall \sigma \in \Sigma.$$

For any $\sigma, \tau \in \Sigma$, $\bar{P}_{\mathrm{sc}}(\sigma, \tau)$ is the conditional probability at stationarity of moving to $\mathcal{X}_\tau$ under $P_{\mathrm{sc}}$, given that the chain is currently in $\mathcal{X}_\sigma$:

$$\bar{P}_{\mathrm{sc}}(\sigma, \tau) = \frac{1}{\pi_{\mathrm{sc}}[\mathcal{X}_\sigma]} \int_{\mathcal{X}_\sigma} \int_{\mathcal{X}_\tau} \pi_{\mathrm{sc}}(x) P_{\mathrm{sc}}(x, dy) \lambda_{\mathrm{sc}}(dx).$$

We will begin by bounding this probability in terms of the probability of moving to $\mathcal{X}_\tau$ under $Q$ (a swap move) and the probability of moving to $\mathcal{X}_\tau$ under $T$ (an update move). For swap moves, let $\bar{Q}$ be the projection matrix of $Q$ with respect to $\{\mathcal{X}_\sigma : \sigma \in \Sigma\}$. Then for $k \in \{0, \ldots, N-1\}$, we have

$$\bar{P}_{\mathrm{sc}}(\sigma, (k, k+1)\sigma) \geq \tfrac{1}{4}\bar{Q}(\sigma, (k, k+1)\sigma) \qquad \forall \sigma$$

where the right-hand side is the conditional probability of swapping $x_{[k]}$ and $x_{[k+1]}$ under $Q$, and then holding twice. Similarly, for update moves we denote by $\bar{T}$ the projection matrix of $T$ with respect to $\{\mathcal{X}_\sigma : \sigma \in \Sigma\}$, and denote $\sigma_{[i,j]} = (\sigma_0, \ldots, \sigma_{i-1}, j, \sigma_{i+1}, \ldots, \sigma_N)$. Then

$$\bar{P}_{\mathrm{sc}}(\sigma, \sigma_{[i,j]}) \geq \tfrac{1}{4}\bar{T}(\sigma, \sigma_{[i,j]}) \qquad \forall i, j.$$

Therefore, the Dirichlet form $\mathcal{E}_{\bar{\mathrm{sc}}}$ of $\bar{P}_{\mathrm{sc}}$ evaluated at $f \in L_2(\pi^*)$ satisfies

$$\mathcal{E}_{\bar{\mathrm{sc}}}(f, f) \geq \tfrac{1}{4}\mathcal{E}_{\bar{Q}}(f, f) + \tfrac{1}{4}\mathcal{E}_{\bar{T}}(f, f).$$

Recalling that $\bar{T}_0$ is the projection matrix of $T_0$ with respect to $\mathcal{A}$, note that

$$4(N+1)\mathcal{E}_{\bar{T}}(f, f)$$

$$= 2(N+1) \sum_{\sigma, \tau \in \Sigma} (f(\sigma) - f(\tau))^2 \pi^*(\sigma) \bar{T}(\sigma, \tau)$$

$$\geq \sum_{\sigma \in \Sigma} \pi^*(\sigma) \sum_{j=1}^J (f(\sigma) - f(\sigma_{[0,j]}))^2 \bar{T}_0(\sigma_0, j)$$

$$= \sum_{\sigma_{1:N}} \left[ \prod_{k=1}^{N} \pi_k[A_{\sigma_k}] \right] \sum_{i=1}^{J} \sum_{j=1}^{J} (f(i, \sigma_{1:N}) - f(j, \sigma_{1:N}))^2 \pi_0[A_i] \bar{T}_0(i, j)$$

$$\geq \sum_{\sigma_{1:N}} \left[ \prod_{k=1}^{N} \pi_k[A_{\sigma_k}] \right] \mathbf{Gap}(\bar{T}_0)$$

$$\times \sum_{i=1}^{J} \sum_{j=1}^{J} (f(i, \sigma_{1:N}) - f(j, \sigma_{1:N}))^2 \pi_0[A_i] \pi_0[A_j]$$

$$= \mathbf{Gap}(\bar{T}_0) \sum_{\sigma \in \Sigma} \pi^*(\sigma) \sum_{j=1}^{J} (f(\sigma) - f(\sigma_{[0,j]}))^2 \pi_0[A_j],$$

where the second inequality is by recognizing the Dirichlet form for $\bar{T}_0$. Therefore,

$$\frac{\mathcal{E}_{\bar{s}c}(f, f)}{\mathbf{Gap}(\bar{T}_0)}$$

$$(16) \qquad \geq \left[ \frac{1}{8} \mathcal{E}_{\bar{Q}}(f, f) + \sum_{\sigma \in \Sigma} \pi^*(\sigma) \sum_{j=1}^{J} (f(\sigma) - f(\sigma_{[0,j]}))^2 \frac{\pi_0[A_j]}{16(N+1)} \right].$$

Now consider a transition kernel $T^*$ constructed as follows: with probability $\frac{1}{2}$ transition according to $\bar{Q}$, or (exclusively) with probability $\frac{1}{2(N+1)}$ draw $\sigma_{[0]}$ according to the distribution $\{\pi_0[A_j] : j = 1, \ldots, J\}$ (i.e., independent samples at the highest temperature); otherwise hold. Note that the Dirichlet form of $T^*$ is exactly four times the right-hand side of (16). Clearly, $T^*$ is also reversible with respect to $\pi^*$, so $\bar{P}_{sc}$ and $T^*$ have the same stationary distribution. Therefore,

$$(17) \qquad \mathbf{Gap}(\bar{P}_{sc}) \geq \frac{\mathbf{Gap}(T^*) \mathbf{Gap}(\bar{T}_0)}{4}.$$

We will now bound $\mathbf{Gap}(T^*)$ by comparison with another $\pi^*$-reversible chain. Define the transition matrix $T^{**}$ which chooses $k$ uniformly from $\{0, \ldots, N\}$ and then draws $\sigma_k$ according to the distribution $\{\pi_k[A_j] : j = 1, \ldots, J\}$. Clearly, $T^{**}$ moves easily among the elements of $\Sigma$, and consequently has a large spectral gap as we will see. By combining (17) with a comparison of $T^*$ to $T^{**}$, we will obtain a lower bound on the spectral gap of $\bar{P}_{sc}$.

Comparison of $T^*$ to $T^{**}$ will be done using Theorem 5.1. To simplify notation, we write $\pi_k(j)$ as shorthand for $\pi_k[A_j]$ for the remainder of this section. Let $j^*$ be the value of $j$ that maximizes $\pi_N(j)$. For each edge $(\sigma, \sigma_{[i,j]})$ in the graph of $T^{**}$, we define a path $\gamma_{\sigma, \sigma_{[i,j]}}$ in $T^*$ with the following 7 stages:

1. Change $\sigma_0$ to $j^*$.
2. Swap that $j^*$ "up" to level $i$.
3. Take the new $\sigma_{i-1}$ (formerly $\sigma_i$) and swap it "down" to level 0.
4. Change the value at level 0 to $j$ (from former $\sigma_i$).
5. Swap the $j$ at level 0 "up" to level $i$.
6. Swap the $j^*$ that is now at level $i-1$ "down" to level 0.
7. Change the value at level 0 to $\sigma_0$ (from $j^*$).

In each path, skip all steps that do not change the state. Using the defined path set, we will obtain an upper bound $c^*$ on the quantity $c$ in Theorem 5.1. Since $T^*$ and $T^{**}$ both have stationary distribution $\pi^*$, Theorem 5.1 then yields $\mathbf{Gap}(T^*) \geq \frac{1}{c^*}\mathbf{Gap}(T^{**})$. To obtain such an upper bound $c^*$, we will use Propositions 6.1 and 6.2.

Proposition 6.1.   *For the above-defined paths,*

$$(18) \qquad \frac{\pi^*(\sigma)T^{**}(\sigma, \sigma_{[i,j]})}{\pi^*(\tau)T^*(\tau, \xi)} \leq 4J\gamma(\mathcal{A})^{-(J+3)}\delta(\mathcal{A})^{-2}$$

*for all $\sigma$, $i$ and $j$, and any edge $(\tau, \xi)$ in $\gamma_{\sigma, \sigma_{[i,j]}}$.*

Proof.   To obtain (18), first note that

$$\pi^*(\sigma)T^{**}(\sigma, \sigma_{[i,j]}) = \frac{\pi_i(j)}{N+1}\left[\prod_{k=0}^{N}\pi_k(\sigma_k)\right]$$

$$= \frac{1}{N+1}\min\{\pi^*(\sigma), \pi^*(\sigma_{[i,j]})\}\max\{\pi_i(\sigma_i), \pi_i(j)\}.$$

For any state $\tau$ in the path $\gamma_{\sigma, \sigma_{[i,j]}}$, we will find a lower bound on $\pi^*(\tau)$ in terms of $\min\{\pi^*(\sigma), \pi^*(\sigma_{[i,j]})\}$. Consider the states in the path $\gamma_{\sigma, \sigma_{[i,j]}}$ up to stage 4 (where $\sigma_i$ is at level 0). We will show that each state $\tau$ satisfies $\pi^*(\tau) \geq \pi^*(\sigma)\gamma(\mathcal{A})^{J+2}J^{-1}$. Then by symmetry, the states from stage 4 to the end of the path satisfy $\pi^*(\tau) \geq \pi^*(\sigma_{[i,j]})\gamma(\mathcal{A})^{J+2}J^{-1}$.

Any state in stages 1 or 2 of the path from $\sigma$ to $\sigma_{[i,j]}$ is of the form $\tau = (\sigma_1, \ldots, \sigma_l, j^*, \sigma_{l+1}, \ldots, \sigma_N)$ for some $l \in \{0, \ldots, i\}$. Therefore,

$$\pi^*(\tau) = \pi^*(\sigma)\left[\prod_{k=1}^{l}\frac{\pi_{k-1}(\sigma_k)}{\pi_k(\sigma_k)}\right]\frac{\pi_l(j^*)}{\pi_0(\sigma_0)}$$

$$= \pi^*(\sigma)\left[\prod_{k=1}^{l}\prod_{m=1}^{J}\left[\mathbf{1}(\sigma_k = m)\frac{\pi_{k-1}(m)}{\pi_k(m)} + \mathbf{1}(\sigma_k \neq m)\right]\right]\frac{\pi_l(j^*)}{\pi_0(\sigma_0)}$$

$$\geq \pi^*(\sigma)\left[\prod_{m=1}^{J}\prod_{k=1}^{N}\min\left\{1, \frac{\pi_{k-1}(m)}{\pi_k(m)}\right\}\right]\frac{\pi_l(j^*)}{\pi_0(\sigma_0)}$$

$$\geq \pi^*(\sigma)\gamma(\mathcal{A})^{J+1}J^{-1},$$

where the last inequality uses the fact that by definition $\pi_N(j^*) \geq J^{-1}$, so $\pi_k(j^*) \geq \gamma(\mathcal{A})J^{-1}$ for all $k$ and

$$\frac{\pi_l(j^*)}{\pi_0(\sigma_0)} \geq \frac{\gamma(\mathcal{A})J^{-1}}{\pi_0(\sigma_0)} \geq \gamma(\mathcal{A})J^{-1}.$$

Now consider the states in stage 3 of the path, the last of which is also the first state in stage 4. Any such state $\tau$ is of the form

$$\tau = (\sigma_1, \ldots, \sigma_l, \sigma_i, \sigma_{l+1}, \ldots, \sigma_{i-1}, j^*, \sigma_{i+1}, \ldots, \sigma_N)$$

for some $l \in \{0, \ldots, i-1\}$. Therefore,

$$\pi^*(\tau) = \pi^*(\sigma)\left[\prod_{k=1}^{l} \frac{\pi_{k-1}(\sigma_k)}{\pi_k(\sigma_k)}\right] \frac{\pi_i(j^*)\pi_l(\sigma_i)}{\pi_0(\sigma_0)\pi_i(\sigma_i)}$$

$$\geq \pi^*(\sigma)\gamma(\mathcal{A})^{J+1}J^{-1}\frac{\pi_l(\sigma_i)}{\pi_i(\sigma_i)}$$

$$\geq \pi^*(\sigma)\gamma(\mathcal{A})^{J+2}J^{-1},$$

where the last step is because $l < i$. Putting the above together, we have that for all $\tau$ in the path from $\sigma$ to $\sigma_{[i,j]}$

$$(19) \qquad \pi^*(\tau) \geq \min\{\pi^*(\sigma), \pi^*(\sigma_{[i,j]})\}\gamma(\mathcal{A})^{J+2}J^{-1}.$$

We use this to obtain (18) as follows: for any edge $(\tau, \xi)$ on the path $\gamma_{\sigma, \sigma_{[i,j]}}$, we have either $\xi = (k, k+1)\tau$ for some $k$, or $\xi = \tau_{[0,m]}$ for some $m$. The probability of proposing the swap $\xi = (k, k+1)\tau$ according to $Q$ is $\frac{1}{2N}$. Recall that for any $l_1, l_2 \in \{1, \ldots, J\}$, $\delta(\mathcal{A})^2$ is a lower bound for the marginal probability at stationarity of accepting a proposed swap between $x_k \in A_{l_1}$ and $x_{k+1} \in A_{l_2}$. Thus, we have $\bar{Q}(\tau, \xi) \geq \delta(\mathcal{A})^2/(2N)$, and so

$$(20) \qquad \begin{aligned} \frac{\pi^*(\sigma)T^{**}(\sigma, \sigma_{[i,j]})}{\pi^*(\tau)T^*(\tau, \xi)} &= \frac{\pi^*(\sigma)\pi_i(j)}{(N+1)\pi^*(\tau)T^*(\tau, \xi)} \\ &\leq \frac{2\pi^*(\sigma)\pi_i(j)}{(N+1)\pi^*(\tau)\bar{Q}(\tau, \xi)} \leq \frac{4\pi^*(\sigma)\pi_i(j)}{\pi^*(\tau)\delta(\mathcal{A})^2} \\ &= \frac{4\min\{\pi^*(\sigma), \pi^*(\sigma_{[i,j]})\}\max\{\pi_i(j), \pi_i(\sigma_i)\}}{\pi^*(\tau)\delta(\mathcal{A})^2} \\ &\leq \frac{4J\max\{\pi_i(j), \pi_i(\sigma_i)\}}{\gamma(\mathcal{A})^{J+2}\delta(\mathcal{A})^2} \leq \frac{4J}{\gamma(\mathcal{A})^{J+2}\delta(\mathcal{A})^2}. \end{aligned}$$

In the case that instead $\xi = \tau_{[0,m]}$ for some $m$, (20) becomes

$$(21) \qquad \frac{\pi^*(\sigma)T^{**}(\sigma, \sigma_{[i,j]})}{\pi^*(\tau)T^*(\tau, \xi)} = \frac{2\pi^*(\sigma)\pi_i(j)}{\pi^*(\tau)\pi_0(m)}$$

and there are three possible cases: the edge $(\tau, \xi)$ could be stage 1, stage 4 or stage 7 of $\gamma_{\sigma, \sigma_{[i,j]}}$. If it is stage 1, then (21) is bounded by

$$\frac{2\pi^*(\sigma)\pi_i(j)}{\pi^*(\sigma)\pi_0(j^*)} \leq \frac{2}{\pi_0(j^*)} \leq \frac{2J}{\gamma(\mathcal{A})}.$$

If the move is stage 4, then (21) is bounded by

$$\frac{2\pi^*(\sigma)\pi_i(j)}{\pi^*(\tau)\pi_0(j)} = \frac{2\min\{\pi^*(\sigma), \pi^*(\sigma_{[i,j]})\}\max\{\pi_i(j), \pi_i(\sigma_i)\}}{\min\{\pi^*(\tau), \pi^*(\xi)\}\max\{\pi_0(j), \pi_0(\sigma_i)\}}$$

$$\leq \frac{2}{\gamma(\mathcal{A})}\frac{\min\{\pi^*(\sigma), \pi^*(\sigma_{[i,j]})\}}{\min\{\pi^*(\tau), \pi^*(\xi)\}} \leq 2J\gamma(\mathcal{A})^{-(J+3)}$$

by (19) and since $\pi_i(j) \leq \frac{\pi_0(j)}{\gamma(\mathcal{A})} \leq \frac{\max\{\pi_0(j), \pi_0(\sigma_i)\}}{\gamma(\mathcal{A})}$ and $\pi_i(\sigma_i) \leq \frac{\pi_0(\sigma_i)}{\gamma(\mathcal{A})} \leq \frac{\max\{\pi_0(j), \pi_0(\sigma_i)\}}{\gamma(\mathcal{A})}$. Finally, if the move is stage 7, then (21) is bounded by

$$\frac{2\pi^*(\sigma)\pi_i(j)}{\pi^*(\sigma_{[i,j]})\pi_0(j^*)} = \frac{2\pi^*(\sigma_{[i,j]})\pi_i(\sigma_i)}{\pi^*(\sigma_{[i,j]})\pi_0(j^*)} \leq \frac{2J}{\gamma(\mathcal{A})}.$$

The result (18) follows for any edge $(\tau, \xi)$ on the path from $\sigma$ to $\sigma_{[i,j]}$. □

PROPOSITION 6.2. *For the above-defined paths,*

$$\sum_{\gamma_{\sigma, \sigma_{[i,j]}} \ni (\tau, \xi)} |\gamma_{\sigma, \sigma_{[i,j]}}| \leq 16(N+1)^2 J^2 \tag{22}$$

*for any edge $(\tau, \xi)$ in the graph of $T^*$.*

PROOF. We will bound the number of paths $\gamma_{\sigma, \sigma_{[i,j]}}$ that go through any edge $(\tau, \xi)$, and the length of any such path.

Consider the set of paths for which the edge is in stage 1 of the path. Then $\tau = \sigma$ and $\xi = \sigma_{[0,j^*]}$, and since $i \in \{0, \dots, N\}$ and $j \in \{1, \dots, J\}$, there are no more than $(N+1)J$ such paths. Similarly, there are no more than $(N+1)J$ paths for which the edge is in stage 4 of the path.

Now consider the set of paths for which the edge is in stage 2 of the path. Then we must have $\tau = (\sigma_1, \dots, \sigma_l, j^*, \sigma_{l+1}, \dots, \sigma_N)$ for some $l \in \{0, \dots, i-1\}$ and $\xi = (l, l+1)\tau$. $\sigma_0$ is unknown but has only $J$ possible values, so with $i$, $j$ unknown there are no more than $(N+1)J^2$ such paths. Similarly, there are no more than $(N+1)J^2$ paths for which the edge is in stage 3 of the path.

If the edge has $\xi = (k, k+1)\tau$ for some $k$, then it can only be in stages 2, 3, 5 or 6 of the path, while if $\xi = \tau_{[0,m]}$ for some $m$ then it can only be in stages 1, 4 or 7. Since the edge can be in at most 4 stages, each with at most $(N+1)J^2$ paths, the total number of paths containing any edge is no more than $4(N+1)J^2$. Each of these paths has length at most $4N+3 < 4(N+1)$, so (22) follows. □

Combining Propositions 6.1 and 6.2, we obtain an upper bound on the constant $c$ in Theorem 5.1:

$$c \leq \frac{2^6(N+1)^2 J^3}{\gamma(\mathcal{A})^{J+3}\delta(\mathcal{A})^2}$$

and recalling that both $T^*$ and $T^{**}$ have stationary distribution $\pi^*$, application of Theorem 5.1 yields

$$\mathbf{Gap}(T^*) \geq \frac{\gamma(\mathcal{A})^{J+3}\delta(\mathcal{A})^2}{2^6(N+1)^2 J^3}\mathbf{Gap}(T^{**}).$$

Now since $T^{**}$ is a product chain whose $\pi_k$-reversible component chains each have spectral gap 1 by definition (1), Theorem 5.3 gives $\mathbf{Gap}(T^{**}) = (N+1)^{-1}$ and we have

$$\mathbf{Gap}(T^*) \geq \frac{\gamma(\mathcal{A})^{J+3}\delta(\mathcal{A})^2}{2^6(N+1)^3 J^3}.$$

Then we obtain the bound for $\mathbf{Gap}(\bar{P}_{\mathrm{sc}})$ from (17):

$$(23) \quad \mathbf{Gap}(\bar{P}_{\mathrm{sc}}) \geq \frac{\mathbf{Gap}(T^*)\,\mathbf{Gap}(\bar{T}_0)}{4} \geq \left(\frac{\gamma(\mathcal{A})^{J+3}\delta(\mathcal{A})^2}{2^8(N+1)^3 J^3}\right)\mathbf{Gap}(\bar{T}_0).$$

Using (14), (15) and (23), then proves Theorem 3.1.

As we have seen, Theorem 3.1 bounds the mixing of the swapping chain in terms of its mixing within each partition element and its mixing among the partition elements. For several multimodal examples, we have given (inverse) temperature specifications which guarantee that the four quantities in the bound are large, and used Theorem 3.1 and Corollary 3.1 to prove rapid mixing of parallel and simulated tempering.

## APPENDIX: PROOF OF THEOREM 5.2

The upper bound in Theorem 5.2 is shown in [14]. The lower bound uses the following results; consider the context of Section 5.2.

THEOREM A.1 [Caracciolo, Pelissetto and Sokal (1992)]. *Let* $\mu_P = \mu_Q$. *Assume that* $P$ *is nonnegative definite and let* $P^{1/2}$ *be its nonnegative square root. Then*

$$\mathbf{Gap}(P^{1/2}QP^{1/2}) \geq \mathbf{Gap}(\bar{P})\min_{j=1,\ldots,J}\mathbf{Gap}(Q|_{A_j}).$$

This result is due to Caracciolo, Pelissetto and Sokal [3], but first appeared in Madras and Randall [14].

Lemma A.1 [Madras and Zheng (2003)].

$$\mathbf{Gap}(P) \geq \frac{1}{n}\mathbf{Gap}(P^n) \qquad \forall n \in \mathbb{N}.$$

Note that although [16] state this result for finite state spaces, their proof extends easily to general spaces.

Lemma A.2 [Madras and Zheng (2003)]. *Assume that $\mu_P = \mu_Q$ and that $P$ is nonnegative definite. Then*

$$\mathbf{Gap}(QPQ) \geq \mathbf{Gap}(P).$$

Now consider the context of Theorem 5.2. The lower bound in Theorem 5.2 follows directly from Theorem A.1 and Lemma A.1:

$$\mathbf{Gap}(P) \geq \frac{1}{2}\mathbf{Gap}(P^2) = \frac{1}{2}\mathbf{Gap}(P^{1/2}PP^{1/2}) \geq \frac{1}{2}\mathbf{Gap}(\bar{P})\min_j \mathbf{Gap}(P|_{A_j}).$$

**Acknowledgments.** We thank the referees for many helpful suggestions.

D. B. Woodard
S. C. Schmidler
Department of Statistical Science
Duke University
Box 90251
Durham, North Carolina 27708-0251
USA
E-mail: dawn@stat.duke.edu
        schmidler@stat.duke.edu

M. Huber
Department of Mathematics
Duke University
Box 90320
Durham, North Carolina 27708-0320
USA
E-mail: mhuber@math.duke.edu